\documentclass[11pt, twoside, a4paper]{article}
\usepackage[english]{babel}
\usepackage{amsmath,amssymb,amsthm,epsfig}
\newtheorem{theorem}{Theorem}

\newtheorem{lemma}{Lemma}
\newtheorem{fact}{Fact}
\newtheorem{definition}{Definition}
\newtheorem{corol}{Corollary}
\theoremstyle{definition}
\newtheorem{remark}{Remark}

\let\s=\sigma

\let\O=\Omega

\newcommand{\R}{\mathbb{R}}
\newcommand{\N}{\mathbb{N}}
\newcommand{\A}{\mathbf{A}}
\newcommand{\Per}[1]{\textrm{Per}_{#1}(\sigma)}
\begin{document}

\title{On the existence of maximizing measures for irreducible countable Markov shifts: a dynamical proof}

\author{Rodrigo Bissacot\thanks{Supported by FAPESP process 2011/16265-8 and CNPq process 454655/2011-8.}\\
\footnotesize{Department of Applied Mathematics, IME-USP, Brazil}\\
\footnotesize{\texttt{rodrigo.bissacot@gmail.com}}
\and
\\
Ricardo dos Santos Freire Jr.\thanks{Supported by FAPESP process 2011/16265-8.} \\
\footnotesize{Department of Mathematics, IME-USP, Brazil}\\
\footnotesize{\texttt{rfreire@usp.br}}
}

\date{\today}
\maketitle
\begin{abstract}
We prove that if  $\Sigma_{\mathbf A}(\N)$ is an irreducible Markov
shift space over $\N$ and \\ $f: \Sigma_{\mathbf A}(\N) \rightarrow \R$
is coercive with bounded variation then there exists a maxi-mizing
probability measure for $f$, whose support lies on a Markov subshift
over a finite alphabet. Furthermore, the support of any maximizing
measure is contained in this same compact subshift. To the best of our knowledge, this is the first proof beyond the finitely primitive case on the  general irreducible non-compact setting. It's also noteworthy that our technique works for the full shift over positive real sequences.
\end{abstract}



\section{Introduction}

Given a dynamical system $T: \Omega \rightarrow \O$ over a space and a real function $f$, the main problem in  \textit{Ergodic Optimization} is to guarantee the existence and characterize the support of the maximizing measures for the system, that is, the invariant Borel probability measures maximizing the operator $\int f d\mu$ over the invariant Borel probabilities for $T$. The survey \cite{Jenkinson} is a good introduction to these problems.

If  $\Omega$ is compact, the existence of the maximizing measures is an immediate consequence of the compactness in the weak*-topology of the set of invariant pro-bability measures. On the other hand, in the non-compact case even the existence is a non-trivial problem. See, for instance, \cite{BG,JMU0,JMU1,JMU2}.

We focus on the case where the space is an {\it irreducible} Markov shift over $\N$ and the dynamics is given by the shift map, that is, $\O=\Sigma_{\A}(\N)$ and $T=\s$. Given $f$ we define $$
 \beta := \sup_{\mu \in \mathcal M_\sigma(\Sigma_{\A}(\N))} \int f \; d\mu \,,$$ 
 and we assume that $f$ is \textit{coercive} in the sense that
\[
\lim_{i \to \infty} \sup f|_{[i]} = - \infty \,,
\]
where $[i]:=\{x \in \Sigma_
{\A}(\N),\,\pi(x)=i \}$ is the cylinder beginning with $i$.  Then, our main result is the following:

\begin{theorem} \label{maintheorem}
Let $\sigma$ be the shift on $\Sigma_{\A}(\N)$ with $\A$ irreducible, $f: \Sigma_{\A}(\N) \to \R$ be a function with bounded variation and coercive. Then, there is a finite set $\mathcal A \subset \N$ such that $\A|_{\mathcal A \times \mathcal A}$ is irreducible and
\[
\beta = \sup_{\mu \in \mathcal M_{\sigma}(\Sigma_{\A}(\mathcal A ))} \int f \ d\mu\,.
\]
Furthermore, if $\nu$ is a maximizing measure, then
\[\text{supp }\nu \subset \mathcal M_{\sigma}(\Sigma_{\A}(\mathcal A ))\,.\]
\end{theorem}

Since $\Sigma_{\A}(\mathcal A)$ is compact, it follows from the first
part of the theorem that there is at least one maximizing measure supported on a subset of $\Sigma_{\A}(\mathcal A)$.

Similar results for {\it finitely primitive}\footnote{The subshift is finitely primitive iff there is $K_0 \in \N$ and a finite sub-alphabet $\mathcal L$ such that any pair of symbols in the alphabet can be connected by a word of exactly $K_0$ of symbols in $\mathcal L$. It's clearly much more stronger than primitive, when you don't require $\mathcal L$ to be finite, which is stronger than irreducible, where there's no uniformity in word length connecting two symbols.} subshifts can be found in \cite{BG,JMU0,JMU1,JMU2,Morris}. In fact, to the best of our knowledge, our result is the first beyond the finitely primitive case, except for the particular case of renewal shifts in \cite{Iommi}.

When $f$ is not coercive, the best known results are \cite{BG,JMU1} that still requires finitely primitive, which follows from the classical oscillation condition. In this case, but in the irreducible context, we're able to prove the following similar result to the ones in \cite{BG,JMU1}:

\begin{theorem} \label{theoremB}
Let $\sigma$ be the shift on $\Sigma_{\A}(\N)$ with $\A$ irreducible, $f: \Sigma_{\A}(\N) \to \R$ be a function with bounded variation and assume there are naturals $I_2 > I_1 > 0$ such that
\[
\sup f|_{[j]} < \beta -\epsilon \quad \forall j \ge I_1\,,
\]
for some $\epsilon>0$ fixed and
\[
\sup f|_{[j]} < C \quad \forall j \ge I_2\,,
\]
where $C$ is a constant (depending on $I_1$) given in \eqref{L}.

Then, there is a finite set $\mathcal A \subset \N$ such that $\A|_{\mathcal A \times \mathcal A}$ is irreducible and
\[
\beta = \sup_{\mu \in \mathcal M_{\sigma}(\Sigma_{\A}(\mathcal A ))} \int f \ d\mu\,.
\]
Furthermore, if $\nu$ is a maximizing measure, then
\[\text{supp }\nu \subset \mathcal M_{\sigma}(\Sigma_{\A}(\mathcal A ))\,.\]
\end{theorem}

Our technique also points out a more natural and elementary approach to the problem of the existence of maximizing measures in the non-compact context. The proofs available up to now pass through the construction of auxiliar functions (normal forms \cite{JMU1} and subactions \cite{BG}) that characterize the support of the maximizing measures or, make use of the thermodynamic formalism \cite{JMU0,Morris}, where more restrictions on the dynamics and the potential $f$ are made. We just use a well-known Parthasarathy's result \cite{Parthasarathy} that says the invariant measures supported on periodic orbits are dense in the ergodic invariant measures for $\s$.

In this way, we reduce our problem to analyzing the ergodic averages of periodic orbits, and the proof is essentially to carry on in details the intuitive idea: since the potencial $f$ decays to $-\infty$ when the symbols grow, we can restrict ourselves to periodic orbits whose symbols are all small.

An important consequence  in these contexts is the
subordination principle, that is a direct application
of the results in \cite{B} or \cite{MorrisSub} after the reduction to the compact
case by our results.

Finally, we remark that our technique can be used in more general contexts, such as the case of the full shift on $\Sigma(\R^+)$.

The paper is organized as follows: in the next section we give the
precise setting and notations to prove the existence part of
the theorems in section 3. In section 4, we finish the proof of the
theorems showing that the support of any maximizing measure must be in
the subshift over the finite alphabet built in the previous section. Finally, in section 5 we point out how our technique works in the case of sequences of positive reals.

\section{Setting and notations}

Let $\N$ be the set of non-negative integers and $\Sigma(\N)$ be the set of sequences of elements in $\N$. Given an infinite matrix $\A: \N \times \N \to \{0,1\}$, we call by $\Sigma_{\A}(\N)$ the subset of $\Sigma(\N)$ of allowable sequences, that is:
\[
\Sigma_{\A}(\N) := \{x \in \Sigma(\N),\, \A(x_i,x_{i+1})=1 \, \forall i\geq 0 \}\,.
\]

Fixed $\lambda \in (0,1)$, we define a metric on $\Sigma_{\A}(\N)$ by $d(x,y) = \lambda^k$, where $k$ is the first coordinate where $x_k \neq y_k$.

Denote by $\pi: \Sigma_{\A}(\N) \to \N$ the projection of the first coordinate, that is $\pi(x) = \pi(x_0 x_1 x_2 \dots) = x_0$.

 We say that $\A$ is \textit{irreducible} when for any $i,j$ in $\N$ there exists a word $w = w_1 \dots w_k$ such that $i w j$ is an allowable word:  $\A(i,w_1)=1$, $\A(w_i,w_{i+1})=1$ for $i=1,\dots,k-1$ and $\A(w_k,j)=1$.

 Our dynamics is given by shift map $\sigma: \Sigma_{\A}(\N) \to \Sigma_{\A}(\N)$ where $(\sigma(x))_i = x_{i+1}$ for all $i\geq 0$ and we denote by $\mathcal M_\sigma (\Sigma_{\A}(\N))$ the set of invariant Borel probability measures for this map. It's clear that $\sigma$ is surjective as $\A$ is irreducible.

Fix a function $f: \Sigma_{\A}(\N) \to \R$ and consider the $j$-th variation of $f$ given by
\[
V_{j}(f) := \sup \{f(x) - f(y) \ , \ \pi(\sigma^i(x))=\pi(\sigma^i(y)) \ \text{for} \  i=0, \dots, j-1\}\,,
\]
and suppose that $f$ has bounded variation, that is
\[
V(f) := \sum_{j=1}^{\infty}V_{j}(f) < \infty\,.
\]

Also, recall that we suppose that $f$ is \textit{coercive} in the sense that
\[
\lim_{i \to \infty} \sup f|_{[i]} = - \infty \,,
\]
where $[i]:=\{x \in \Sigma_
{\A}(\N),\,\pi(x)=i \}$ is the cylinder beginning with $i$.

Since $f$ is coercive and has bounded variation, it's easy to see that $f$ is continuous and bounded above, which implies that $\beta$ as defined in the introduction is, in fact, well defined.

Our existence problem is to show there is a finite alphabet $\mathcal A \subset \N$ and a maxi\-mizing measure for $f$, that is, an invariant probability measure $\nu\in \mathcal M_\sigma(\Sigma_{\A}(\N))$ such that
\[
\beta = \int f \; d\nu \,,
\]
where $\nu$ is supported on $\Sigma_{\A}(\mathcal A)$, the set of allowable sequences of symbols in $\mathcal A$.


\section{Proof of the existence results}
Let $\mathcal M_{\sigma-Per}(\N)$ be the set of periodic invariant probability measures, that is, the invariant probability measures that are supported on a periodic orbit of $\sigma$. This set is extremely important since we can reduce the problem into the study of periodic orbits through the following lemma.

\begin{lemma} \label{lemmabeta}
\[
\beta = \sup_{\mu \in \mathcal M_{\sigma-Per}(\Sigma_{\A}(\N))} \int f \; d\mu \,.
\]
\end{lemma}
\proof
The Ergodic Decomposition theorem implies that
\[
\beta = \sup_{\mu \in \mathcal M_{\sigma-erg}(\Sigma_{\A}(\N))} \int f \; d\mu \,,
\]
where $\mathcal M_{\sigma-erg}(\Sigma_{\A}(\N))$ is the set of ergodic invariant probability measures.



By \cite{Parthasarathy}
the periodic invariant probability measures are dense in  $\mathcal M_{\sigma-erg}(\Sigma_{A}(\N))$ and we're done.
\qed

We denote the set of $n$-periodic points of $\sigma$ by $\Per{n}$ and the set of all $\sigma$-periodic orbits is $\text{Per}(\sigma):=\bigcup_{n \ge 1} \Per{n} $.

\begin{definition}
Let $S_mf(x) := \displaystyle\sum_{i=0}^{m-1}f(\sigma^i(x))$, we use the following notation:
\begin{enumerate}
\item[i)] for any $x \in \Sigma_\A(\N)$ let $\beta_m(x):= \displaystyle \frac{1}{m}S_mf(x)$;
\item[ii)] for any $x \in \Sigma_\A(\N)$ denote by $\beta(x):=\displaystyle
  \lim_{m\to \infty} \beta_m(x) $ whenever the limit exists. Notice that if $x \in \Per{n}$ then the ergodic average of $x$ is
  $\beta(x) = \beta_n(x)$;
\item[iii)] we say that $ x \in \Sigma_\A(\N)$ starts in $i$ when $i$ is
  the smallest natural that appear in the coordinates of $x$. In
  particular, all symbols of $x$ are greater or equal to $i$ and if $x \in \Per{n}$ we have
\[
i = \min_{0 \le j \le n-1}\{\pi(\sigma^j(x)) \}\,;
\]
\item[iv)] given a pair $i, j$ in $\N$ we say that a word $w = w_1 \dots w_k$ connects $i$ to $j$ when $i w j$ is an allowable word: $\A(i,w_1)=1$, $\A(w_i,w_{i+1})=1$ for $i=1,\dots,k-1$ and $\A(w_k,j)=1$.
\end{enumerate}
\end{definition}

Since $\int f \ d\mu = \frac{1}{n}\sum_{i=0}^{n-1}f(\s^i(x))$ when $\mu$ is a periodic invariant probability supported on the orbit of $x \in \Per{n}$, it's clear from lemma \ref{lemmabeta} that $$\beta = \sup_{x \in
  \text{Per}(\sigma)}\beta(x)\,.$$ Then, our problem is reduced into showing the existence of a finite alphabet $\mathcal A \subset \N$ such that
\begin{equation} \label{paparazi}
\beta = \sup_{x \in \text{Per}(\sigma) \cap \Sigma_{\A}(\mathcal A)}\beta(x) \,.
\end{equation}

Now we make a first cut on the symbols. The following lemma, together with lemma \ref{lemmabeta}, implies that we don't have to care about periodic orbits whose symbols are all too large.

\begin{lemma} \label{lemmaI1}
Given $\epsilon >0$, there is $I_1 \in \N$ such that if $x$ starts in
$i \ge I_1$ then $\beta_m(x) < \beta -\epsilon$ for any $m \in \N$. In
particular, if $x \in \Per{n}$ we have $\beta(x) < \beta -\epsilon$.
\end{lemma}
\proof
Since $f$ is coercive, there is $I_1 \in \N$ such that
\[
\sup f|_{[j]} < \beta -\epsilon \quad \quad \text{for all } j \ge I_1 \,.
\]

We have that
\[
\beta_m(x) = \frac{1}{m}S_mf(x) = \frac{1}{m}\sum_{j=0}^{m-1}f(\sigma^j(x)) \le \frac{1}{m} \sum_{j=0}^{m-1} \sup f|_{[\pi(\sigma^j(x))]} \,,
\]
and since $\pi(\sigma^j(x)) \ge i \ge I_1$ for all $j=0, \dots, m-1$, we get
\eject
\[
\beta_m(x) \le \frac{1}{m} \sum_{j=0}^{m-1} \sup f|_{[\pi(\sigma^j(x))]} < \beta -\epsilon \,.
\] \qed

Let us fix $\epsilon > 0$. If we consider the alphabet $\mathcal I_1 := \{0, 1, \dots, I_1-1\}$ we still have a problem that maybe there are no allowable sequences only with such symbols and, besides, the shift does not need to be irreducible when restrict to such sequences. So we complete $\mathcal I_1$ to a finite alphabet $\mathcal A_1$ in the following manner.

We choose, for each pair $i,j$  in $\mathcal I_1$, one word $w=w(i,j)$ connecting $i$ to $j$. Notice there is such a word since $\A$ is irreducible. We denote by $P_0$ the length of the longest of such connecting words. Let $\mathcal C_1$ be the set of symbols that appear in at least one of these connecting words and then consider $\mathcal A_1:= \mathcal I_1 \cup \mathcal C_1$. Since each connecting word has at most a finite number of symbols, and we have chosen $I_1^2$ words, we have that $\mathcal A_1$ is finite.

It's clear that any pair of symbols in $\mathcal A_1$ can be connected using only symbols in $\mathcal A_1$. This means that $\A$ restricted to $\mathcal A_1$ is irreducible.

Therefore $\Sigma_{\A}(\mathcal A_1):= \{x \in \Sigma_{\A}(\N), \, \pi(\sigma^i(x)) \in \mathcal A_1 \ \forall i\ge 0 \}$ is a compact invariant subspace of $\Sigma_{\A}(\N)$.

Now we can make a second cut on the alphabet and show it's enough. In fact, since $f$ is coercive, there is $I_2 \ge I_1$ such that

\begin{equation} \label{L}
\begin{split}
\sup f|_{[j]} <  C := \min \left\{ C_1, C_2 \right\} \quad \forall j\ge I_2 \,,
\end{split}
\end{equation}
where
\begin{align*}
C_1 := & - \left( P_0|\min f|_{\Sigma_{\A}(\mathcal A_1)}| + (P_0-1) |\beta| + 2V(f) \right) \,, \\
C_2 := & \beta-\epsilon - V(f) \,.
\end{align*}

Below, in the proof of our key lemma \ref{lemmafinale}, we create a new periodic orbit with smaller symbols by connecting two symbols in $\mathcal A_1$ appearing on the orbit. First, in the case when we have chosen a non empty connecting word we need the estimate given by $C_1$. Otherwise, we have an empty connecting word, that is, we can connect both symbols directly, and we need the estimate given by $C_2$. 

Then we can complete $\mathcal I_2 = \{0, 1, \dots, I_2-1\}$ into a
finite alphabet $\mathcal A_2$ in the same way we did with $\mathcal A_1$, with the same dynamical properties. It's also clear that we can take $\mathcal A_2$ such that $\mathcal A_1 \subset \mathcal A_2$.

We need some control over the ergodic average on parts of a given orbit. For that purpose, the following definition is convenient:

\begin{definition}
Let $x \in \Sigma_{\A}(\N)$ and $w = x_\ell \dots x_{\ell+m}$ be a word appearing on $x$. Then:
\begin{enumerate}
\item the ergodic average of the word $w = x_\ell \dots x_{\ell+m}$ on the orbit $x$ is
\[ \kappa(\ell, m | x) = \kappa(w | x) := \frac{1}{m+1} \sum_{j=0}^{m}f(\sigma^{\ell+j}(x)) \,;\]
\item if $r<m$ we define
\[ \kappa_r(\ell, m | x) = \kappa_r(w | x) := \frac{1}{r+2} \left( f(\sigma^{\ell+m}(x)) + \sum_{j=0}^{r}f(\sigma^{\ell+j}(x))\right) \,.\]
\end{enumerate}
\end{definition}

The following facts shows the relation between the previous definition and the ergodic average of a periodic orbit.

\begin{fact} \label{factk}
Let $x \in \Sigma_{\A}(\N)$ be a periodic orbit for $\sigma$ such that $\beta(x) \ge \beta-\epsilon$ and $x \notin \Sigma_{\A}(\mathcal A_2)$. Then, there is at least one word $x_{\ell} \dots x_{\ell+m}$ appearing in $x$ such that
\begin{enumerate}
\item $\kappa(\ell, m | x) \ge \beta(x)$;
\item $x_{\ell} < I_1$, $x_{\ell+m}\ge I_2$; and
\item $x_{\ell+j} < I_2$ for all $j \in \{0, \dots, m-1\}$.
\end{enumerate}
\end{fact}
\proof
Since $\beta(x) \ge \beta-\epsilon $, lemma \ref{lemmabeta} implies that $x$ starts in $i < I_1$.

And because $x \notin \Sigma_{\A}(\mathcal A_2)$, we have that there is at least one symbol greater or equal to $I_2$ appearing on $x$, as by construction we have that $\mathcal I_2 \subset \mathcal A_2$.

This shows that there is at least one word appearing in $x$ satisfying both properties 2 and 3. For each such word, we may take it to be the longest one satisfying such properties, and in this sense let us call it a maximal word.

Since $x$ is periodic, there is at most a finite number of such maximal words appearing in $x$. Also, if a symbol on $x$ is not on any of these maximal words, it must
be greater or equal to $I_1$, otherwise it would be possible to
extend a maximal word, which is absurd.

We can write a period of $x$ as a concatenation of maxi\-mal and
non maxi\-mal words, that is, $w_0 \dots w_k$ represents a period of $x$
and each word $w_j$ for $j \in \{0, \dots, k \}$ is either maximal or
has only symbols greater or equal to $I_1$. Let $\ell_j$ be the length
of the word $w_j$ and we get
\[
\beta(x) = \frac{1}{n} \sum_{j=0}^{k} \ell_j \kappa(w_j|x) \,.
\]

Let $\bar\ell$ be a word $w_{\bar\ell}$ such that $\kappa(w_{\bar\ell}|x) = \max_{j \in \{0, \dots, k \}} \{
\kappa(w_j|x)\}$. Then

\[
\beta(x) \le \frac{1}{n} \sum_{j=0}^{k} \ell_j \kappa(w_{\bar\ell}|x) = \kappa(w_{\bar\ell}|x) \,.
\]
As in lemma \ref{lemmaI1}, if $w_{\bar\ell}$ is not
one of the maximal words, then $\kappa(w_{\bar\ell}|x) < \beta -\epsilon $. Since $\beta(x) \ge \beta - \epsilon$, $w_{\bar\ell}$ must be one of the maximal
words and our claim follows
taking $x_\ell \dots x_{\ell+m} := w_{\bar\ell}$.
\qed

\begin{fact}  \label{factkr}
Let $x_{\ell} \dots x_{\ell+m}$ be the word given by fact 1 and $r<m$ be the greatest integer such that $x_{\ell+r} \in \mathcal I_1$. Then
\[\kappa(\ell,m|x) \le \kappa_r(\ell,m|x) \,. \]
\end{fact}
\proof
In fact, we have by definition that
\[
\kappa(\ell,m|x)  =  \frac{1}{m+1}\left( (r+2)\kappa_r(\ell,m|x) + \sum_{j=r+1}^{m-1}f(\sigma^{\ell+j}(x))\right)
\]
and since $\pi(\sigma^{\ell+j}(x)) \ge I_1$ for $j\ge r+1$,
from the same argument of lemma \ref{lemmaI1} it follows that
\[
 \kappa(\ell,m|x) \le  \frac{1}{m+1}\left( (r+2)\kappa_r(\ell,m|x) + (m-r-1) (\beta -\epsilon)  \right)\,,
\]
and recall from fact \ref{factk} that $\kappa(\ell,m|x) \ge \beta(x) \ge \beta -\epsilon $ and so
\[
 \kappa(\ell,m|x) \le \frac{r+2}{m+1} \ \kappa_r(\ell,m|x) + \frac{m-r-1}{m+1} \ \kappa(\ell,m|x) \,.
\]

Now, reordering the last expression
\[
\frac{r+2}{m+1} \ \kappa(\ell,m|x) \le \frac{r+2}{m+1} \ \kappa_r(\ell,m|x) \,,
\]
and the result follows.
\qed

Let $\delta := \min \{ C_1, C_2 \} - \sup f|_{\cup_{j\ge I_2}[j]} > 0$, where the constants are from
\eqref{L}. The following lemma is the key to complete the proof of
theorem \ref{maintheorem}.

\begin{lemma} \label{lemmafinale}
Let $x \in \Sigma_{\A}(\N)$ be any periodic orbit for $\sigma$ such
that $x \notin \Sigma_{\A}(\mathcal A_2)$ and $\beta(x) \ge
\beta-\epsilon$. Then, there is a periodic orbit  $z \in
\Sigma_{\A}(\mathcal A_2)$ such that $\beta(z) > \beta(x)$.
\end{lemma}
\proof
Consider $x_{\ell} \dots x_{\ell+m}$ the word given by fact \ref{factk}, and $r<m$ the greatest integer such that $x_{\ell+r} \in \mathcal I_1$.

Now, take $z = (\overline{ x_{\ell} \dots x_{\ell+r} w })$, that is, the orbit made by repetition of the word $x_{\ell} \dots x_{\ell+r} w$, where $w$ is the word of size $q$ connecting $x_{\ell+r}$ to $x_{\ell}$ made of symbols in $\mathcal A_1$, chosen in the definition of $\mathcal A_1$.

Notice that both $x_{\ell}$ and $x_{\ell+r}$ are in $\mathcal I_1$ by
facts \ref{factk} and \ref{factkr}, but $x_{\ell+r+1}$ may not be in
$\mathcal A_1$, and it's important for our estimates bellow that we
use only connecting symbols in $\mathcal A_1$.

By facts \ref{factk} and \ref{factkr}, we know that $\kappa_r(\ell,m|x) \ge \kappa(\ell,m | x) \ge
\beta(x)$. So, we're left to show that $\beta(z) -\delta_1 \ge
\kappa_r(\ell,m|x)$ for some $\delta_1>0$.

In fact, we have
\[
\kappa_r(\ell,m|x) = \frac{1}{r+2} \left( f(\sigma^{\ell+m}(x)) + \sum_{j=0}^{r}f(\sigma^{\ell+j}(x))\right) \,,
\]
and since $f$ has bounded variation and $x_{\ell+j}=z_j$ for all $j \in \{0, \dots, r \}$, we have that
\[
\sum_{j=0}^{r}f(\sigma^{\ell+j}(x)) \le \sum_{j=0}^{r}f(\sigma^{j}(z)) + V(f)\,,
\]
so we get
\begin{equation} \label{kr}
\kappa_r(\ell,m|x) \le \frac{1}{r+2} \left( f(\sigma^{\ell+m}(x)) + \sum_{j=0}^{r}f(\sigma^{j}(z)) + V(f)\right) \,.
\end{equation}

Recall that $\pi(\sigma^{\ell+m}(x)) \ge I_2$ and from
\eqref{L} and the definition of $\delta$ we have that

\begin{equation} \label{L2}
f(\sigma^{\ell+m}(x)) \le \sup f|_{\cup_{j\ge I_2}[j]} = \min \{C_1, C_2 \} - \delta \le C_i -\delta\,,
\end{equation}
for $i=1,2$. \\

We have 2 cases to consider: $q \ge 1$ and $q = 0$. \\


First, assume that $q \ge 1$ and recall from
\eqref{L2} and $C_1$ in \eqref{L} that
\begin{equation*}
f(\sigma^{\ell+m}(x)) \le -P_0 |\min f|_{\Sigma_{\A}(\mathcal A_1)}| -
2V(f) + (1-P_0)|\beta| - \delta \,.
\end{equation*}

Notice that
\begin{eqnarray*}
\sum_{j=1}^{q}f(\sigma^{r+j}(z)) & \ge & \sum_{j=0}^{q-1}f(\sigma^{j}(z_{\text{aux}})) -V(f) \\
& \ge &  q \min f|_{\Sigma_{\A}(\mathcal A_1)} - V(f) \\
& \ge & -P_0 |\min f|_{\Sigma_{\A}(\mathcal A_1)}| - V(f) \,,
\end{eqnarray*}
where $z_{\text{aux}}$ is any point in $\Sigma_{\A}(\mathcal A_1)$ starting by the word $w$. For example, we can take $z_{\text{aux}}$ a periodic point, connecting $w_q$ to $w_1$ just like we did to obtain $z$.

In this way, we get in the previous inequality
\begin{equation*}
f(\sigma^{\ell+m}(x)) \le  \sum_{j=1}^{q}f(\sigma^{r+j}(z))-V(f) + (1-P_0)|\beta| - \delta \,.
\end{equation*}

Applying this to \eqref{kr} we have

\begin{eqnarray*}
\kappa_r(\ell,m|x) & \le &  \frac{1}{r+2} \left(
  \sum_{j=0}^{r} f(\sigma^{j}(z)) + \sum_{j=1}^{q} f(\sigma^{r+j}(z)) + (1-P_0)|\beta| - \delta \right) \\
& = & \frac{r+1+q}{r+2} \beta(z) +
\frac{1-P_0}{r+2} |\beta| - \frac{\delta}{r+2} \\
& = & \beta(z) + \frac{1}{r+2} \left( (q - 1)\beta(z) + (1-P_0)|\beta| - \delta  \right) \,,
\end{eqnarray*}
and since $\beta(z) \le \beta \le |\beta|$, we have that
\[
(q - 1)\beta(z) \le (q - 1)|\beta| \le (P_0 -1) |\beta|\,,
\]
implying that $(q - 1)\beta(z) + (1-P_0)|\beta|
\le 0$. Therefore  $$\kappa_r(\ell,m|x) \le \beta(z) -
\frac{\delta}{r+2} := \beta(z) - \delta_1\,,$$ as we wanted.


Finally, assume that $q = 0$. That means $z= (\overline{
  x_{\ell} \dots x_{\ell+r}})$.

From \eqref{L2} and $C_2$ in \eqref{L} we have that
\[
f(\sigma^{\ell+m}(x)) \le \beta - \epsilon-V(f) - \delta\,,
\]
and from \eqref{kr} we get
\begin{eqnarray*}
\kappa_r(\ell,m|x) & \le  & \frac{1}{r+2} \left( f(\sigma^{\ell+m}(x)) + \sum_{j=0}^{r}f(\sigma^{j}(z)) + V(f)\right) \\
& \le & \frac{1}{r+2} \left( \sum_{j=0}^{r}f(\sigma^{j}(z))  + \beta - \epsilon - \delta\right) \\
& = & \frac{1}{r+2} \left( (r+1) \beta(z) + \beta - \epsilon - \delta\right) \\
& = & \beta(z) + \frac{\beta - \epsilon - \beta(z)}{r+2} - \frac{\delta}{r+2}   \,,
\end{eqnarray*}
and since by facts \ref{factk} and \ref{factkr} we have
$\kappa_r(\ell,m|x)  \ge \beta - \epsilon$, the
last inequality also implies that
\[
\beta(z) + \frac{\beta - \epsilon - \beta(z)}{r+2} \ge \beta -
\epsilon \,,
\]
from which we have that $\beta - \epsilon - \beta(z) \le 0 $, and so $$\kappa_r(\ell,m|x) \le \beta(z) -
\frac{\delta}{r+2} = \beta(z) - \delta_1\,, $$ as desired.
\qed


\begin{remark}
It's important to realize that in lemma \ref{lemmafinale} we've proved that exchanging $x_{\ell+m} \ge I_2$ for $w$, we have increased by at least $\delta>0$ in the ergodic sums. That is, $S_{m+1}f(\s^{\ell}(x)) + \delta \le S_{p}f(z)$, where $p$ is the period of $z$. This will be important in the next section.
\end{remark}

Now we're able to complete the proofs of the existence of a maximizing
measure.

\proof[Proof of the existence in Theorem \ref{maintheorem}]

Recall that our problem is reduced into proving \eqref{paparazi}, and that lemma \ref{lemmabeta} implies that
\[
\beta = \sup_{x \in \text{Per}(\sigma)} \beta(x) \,.
\]

Let $x^n \in \text{Per}(\s)$ for all $n$ be a sequence of periodic orbits such that $\beta(x^n) \to \beta$ as $n\to \infty$ and, so, we can assume $\beta(x^n) \ge \beta - \epsilon$.

We take $\mathcal A = \mathcal A_2$ as defined before, and then lemma
\ref{lemmafinale} shows that, for each $n$ there is a periodic point
$z^n \in \Sigma_{\A}(\mathcal A)$ such that $\beta(z^n) \ge \beta(x^n)$.\footnote{Notice that $z^n$ here may be taken
as $x^n$ if $x^n \in \Sigma_{\A}(\mathcal A)$.} Therefore, as $\beta(x^n) \to \beta$ as $n\to \infty$, so does $\beta(z^n)$, and we're done. \qed

\proof[Proof of the existence in Theorem \ref{theoremB}]

The theorem follows in the same way. In fact, in the proof of theorem \ref{maintheorem} we only use the fact that $f$ is coercive to guarantee the existence of $I_1$ and $I_2$ satisfying the hypothesis given, and to guarantee that $\delta > 0$. In this case, it's enough to consider $\delta := \min \{ C_1, C_2 \} - \sup f|_{[x_{\ell+m}]} > 0 $.
\qed


\section{Proof that $\text{supp } \nu \subset \Sigma_\A(\mathcal A)$
  for any $\nu$ maximal}
Let us keep the same notation from the previous section, in particular recall that $\mathcal A = \mathcal A_2$. The proof for
theorems \ref{maintheorem} and \ref{theoremB} are similar, so we make
no distinction here.

We know from the previous section that there is at least one maximal
measure whose support is in $\Sigma_\A(\mathcal A_2)$. Besides, from
lemma \ref{lemmafinale}, we also know that there is no periodic
maximal measure whose support is not contained in $\Sigma_\A(\mathcal
A_2)$.

Now consider $\nu$ a non periodic maximal measure for $f$ and by contradiction suppose that $\text{supp }
\nu \not\subset \Sigma_\A(\mathcal A_2)$. By the Ergodic Decomposition
theorem we can suppose $\nu$ is ergodic.

The key step now is to build an invariant periodic measure using a
generic point on $\text{supp } \nu$, whose
ergodic average is strictly greater than $\beta$, which is absurd and
proves our result. It's convenient to consider
$\beta = 0 $ here.\footnote{This is easily done by considering $f -
  \beta$.}

Let $x = (x_0 x_1 x_2 \dots) \in \text{supp } \nu$ be a generic point such that $\beta_m(x)
\to 0$ as $m \to \infty$ and $x$ is recurrent. Because
of lemma \ref{lemmaI1}, we can assume without loss of generality that $x_0 < I_1$.

Since $\text{supp }
\nu \not\subset \Sigma_\A(\mathcal A_2)$, there
is a symbol $I \ge I_2$ appearing in the expression of $x$.

We want to modify $x$ into a new point $z$ such that $z$ is periodic
and $\beta(z) > 0$. Since this periodic orbit induces an invariant
periodic measure $\mu$ that has $\int f \ d\mu = \beta(z) > 0$, this
gives a contradiction with the fact that we took $\nu$ a maximizing
measure, and we're done. Notice that there's no need for $z$ to be in
$ \Sigma_\A(\mathcal A_2) $ for this to work.

Let $i$ be the smallest integer such that $x_i=I$, and consider $b <
i$ the greatest integer such that $x_b < I_1$ and $a > i$ the smallest
integer such that $x_a < I_1$.

In this way, we find a word $w$ beginning with $x_b$ and ending in
$x_a$ such that between them there are only symbols greater or equal
to $I_1$ and at least one symbol equal to $I \ge I_2$. We aim at
exchanging $w$ for $\widetilde w$, which is another word beginning in $x_b$
and ending in $x_a$ but between them we put, if necessary, a
connecting word $y$ made of symbols in $\Sigma_\A(\mathcal A_1)$. That
is, $\widetilde w = x_b y x_a$ (but maybe $\widetilde w = x_b x_a$ ).

Let $m_1$ be the size of the prefix of $x$ finishing precisely after
the first appearance of $w$ in the expression of $x$. Exchanging $w$ for
$\widetilde w$, and using the same calculations\footnote{Notice that,
  in fact, the calculations are a little bit easier here since we're
  only looking for the ergodic sums and not the averages. See remark 1.} in the proof of lemma
\ref{lemmafinale}, we get $x^1$ with a modified new prefix of size
$\widetilde m_1$ such that $S_{m_1}f(x) + \delta \le S_{\widetilde m_1}f(x^1)$.

We can repeat this process with the next appearances of $w$ in the
expression of $x$, and after $k$ exchanges, we get a new point $x^k$
such that $$ S_{m_k}f(x) + k \delta \le S_{\widetilde m_k}f(x^k). $$
Then, let $N$ be an integer such that $(N-1) \delta \ge P_0|\min
f|_{\Sigma_\A(\mathcal A_1)}| + 2 V(f)$. Therefore, we have
\begin{equation}\label{S_m_N}
  S_{\widetilde m_N}f(x^N) \ge S_{m_N}f(x) +
  \delta + P_0|\min f|_{\Sigma_\A(\mathcal A_1)}| + 2 V(f)\,.
\end{equation}

Applying Atkinson's lemma \cite{ATK}, we get $\ell \ge m_N$ such that
$|S_\ell f(x)| \le \frac{\delta}{2}$. Fix $m \le \ell$ as the greatest
integer such that $x_m < I_1$. It's clear that $m \ge m_N$. Also,
since $x_{m+1}, \dots, x_\ell$ are greater or equal to $I_1$, from the
choice of $I_1$ on lemma \ref{lemmaI1}, we have that $S_m f(x) \ge
S_\ell f(x)$ and so we get
\begin{equation}\label{atkinson}
 S_m f(x) \ge -\frac{\delta}{2} \,.
\end{equation}

Notice that $S_m f(x) = S_{m_N} f(x) + S_{m-m_N } f( \sigma^{m_N}(x)
)$. Let $\widetilde m$ is the position of $x_m$ in $x^N$ after the $N$
exchanges we made, and since by the definition of $x^N$ we have $\sigma^{m_N}(x) =
\sigma^{\widetilde m_N}(x^N)$, we also have that $\widetilde m -
\widetilde m_N = m - m_N$. Therefore
\[S_{\widetilde m} f(x^N) =
S_{\widetilde m_N} f(x^N) + S_{m-m_N } f( \sigma^{m_N}(x)
)
\]
and using \eqref{S_m_N} we get
\begin{equation}\label{closing}
  S_{\widetilde m} f(x^N) \ge S_{m}f(x) + \delta + P_0|\min
  f|_{\Sigma_\A(\mathcal A_1)}| + 2 V(f)\,.
\end{equation}

Now, let $u$ be a word in $\Sigma_\A(\mathcal A_1)$ of size $q$ connecting
$x_m$(=$x_{\widetilde m}$) to $x_0$, and let $z \in \Sigma_\A(\N)$ be
the periodic point given by the repetition of the word $x^N_0 x^N_1
\dots x^N_{\widetilde m -1} u$, that is, the prefix of size
$\widetilde m$ of $x^N$ concatenated with $u$. Let $p = \widetilde m +
q$ be the period of $z$.

Similarly to the proof of lemma \ref{lemmafinale}, we have that
\begin{eqnarray*}
  S_p f(z) & = & S_{\widetilde m} f (z) + S_q f (\sigma^{\widetilde
    m}(z)) \\
& \ge & S_{\widetilde m} f (x^N) - V(f) + S_q f (\sigma^{\widetilde
    m}(z))\\
& \ge & S_{\widetilde m} f (x^N) - 2V(f)  -P_0 |\min
f_{\Sigma_\A(\mathcal A_1)} | \,,
\end{eqnarray*}
and by \eqref{closing} we get
\[
S_p f(z) \ge S_{m}f(x) + \delta \,,
\]
so that by \eqref{atkinson} we have
\[
S_p f(z) \ge \frac{\delta}{2} \,,
\]
and we're finally done, since
\[
\beta(z) = \frac{1}{p}S_p f(z) \ge \frac{\delta}{2p} > 0 \,.
\]

\section{The case of $\Sigma(\R^+)$}
In the case of the full shift $\sigma$ on $\Sigma(\R^+):= {\R^+}^{\N} $, the sequences of positive reals, where the shift is the same as before, and all sequences are allowable, the previous technique works. In fact the proof of the theorem is easier, since the proof of lemma \ref{lemmafinale} is restricted to the case when we don't need any further symbols to create the orbit $z$.

In particular, we only need to consider the second constant in \eqref{L}.

Finally, we notice that in this case the fact that $f$ is coercive and has bounded variation does not imply that $f$ is bounded above, so we have to make this hypothesis to guarantee the existence of $\beta$.

In this way, we get the following theorem, that is an analogous to corollary 6.2 in \cite{JMU1}:

\begin{theorem} \label{realcase}
Let $\sigma$ be the full shift on $\Sigma(\R^+)$ and $f: \Sigma(\R^+) \to \R$ be a bounded above function with bounded variation and assume there are real numbers $I_2>I_1>0$ such that
\[
\sup f|_{[j]} < \beta-\epsilon \quad \forall \ j \ge I_1 \,,
\]
for some $\epsilon > 0$ fixed and
\[
\sup f|_{[j]} < \min f|_{\Sigma([0,I_1])} - V(f) \quad \forall \ j \ge I_2 \,.
\]

Then, we have that
\[
\beta = \sup_{\mu \in \mathcal M_{\sigma}(\Sigma([0,I_2]))} \int f \ d\mu\,.
\]
Furthermore, if $\nu$ is a maximizing measure, then
\[\text{supp }\nu \subset \mathcal M_{\sigma}(\Sigma([0,I_2]))\,.\]
\end{theorem}

In the case when $f$ is coercive, we have the following:

\begin{corol} \label{realcoercivecase}
Let $\sigma$ be the full shift on $\Sigma(\R^+)$ and $f: \Sigma(\R^+) \to \R$ be a bounded above function with bounded variation and coercive. Then, there is $I > 0$ such that
\[
\beta = \sup_{\mu \in \mathcal M_{\sigma}(\Sigma([0,I]))} \int f \ d\mu\,.
\]
Furthermore, if $\nu$ is a maximizing measure, then
\[\text{supp }\nu \subset \mathcal M_{\sigma}(\Sigma([0,I]))\,.\]
\end{corol}


\section*{Acknowledgments}

The authors are very grateful to the Professor F\'{a}bio A. Tal for all the helpful comments and suggestions, especially for pointing out Atkinson's lemma to us, and to Professor Rafael Rig\~{a}o Souza for the careful reading and feedback.

\end{document}